\documentclass[10pt,journal]{IEEEtran}
\usepackage{graphicx}          
\usepackage{amsmath}           
\usepackage{amssymb}
\usepackage{algorithm}
\usepackage{algpseudocode}
\usepackage[latin1]{inputenc}   

\usepackage{amsthm}

\theoremstyle{definition}

\DeclareMathOperator{\E}{E}
\DeclareMathOperator*{\argmin}{arg\,min}

\begin{document}

\title{Alternating Least-Squares for \\ Low-Rank Matrix Reconstruction}
\author{Dave Zachariah, Martin Sundin, Magnus Jansson and Saikat Chatterjee\thanks{The authors are with the ACCESS Linnaeus Centre, KTH Royal Institute of
Technology, Stockholm. E-mail:
$\{$dave.zachariah,magnus.jansson$\}$@ee.kth.se, masundi@kth.se and
saikatchatt@gmail.com}}

\maketitle
\begin{abstract}
For reconstruction of low-rank matrices from undersampled measurements, we develop an iterative algorithm based on  least-squares estimation.
While the algorithm can be used for any low-rank matrix, it is also capable of exploiting a-priori knowledge of matrix structure. In particular, we consider linearly structured matrices, such as Hankel and Toeplitz, as well as positive semidefinite matrices.
The performance of the algorithm, referred to as alternating least-squares (ALS), is evaluated by simulations and compared to the Cramér-Rao bounds.
\end{abstract}


\section{Introduction}

Low-rank matrices appear in various areas of signal processing and system identification \cite{Werner&Jansson2006}. In recent times the problem of reconstructing such matrices from a vector of linear measurements embedded in noise has attracted a lot of attention \cite{Candes&Plan2010}. In this scenario the dimension of the measurement vector is lower than the number of matrix elements, and hence the problem consists of an underdetermined set of linear equations. This problem has a myriad of applications, for example, spectral imaging \cite{SignerettoEtAl2011}, wireless sensor networks \cite{ChengEtAl2010}, video error concealment \cite{DaoEtAl2010}.

In the gamut of designing low-rank matrix reconstruction algorithms, most of the existing research work deals with a specific setup known as `matrix completion' where the measurement vector consists of a subset of elements of the underlying matrix. The algorithms can be separated into three broad categories: Convex-relaxation based \cite{CaiEtAl2008,Toh&Yun2010}, minimization on Grassmannian manifold of subspaces \cite{BalzanoEtAl2010,DaiEtAl2010},
least-squares matrix fitting \cite{WenEtAl2010,Tang&Nehorai2011}.
While we note that a substantial effort is paid to the matrix completion problem, far less effort is devoted to a general setup consisting of an underdetermined system of linear equations. There are few attempts, such as \cite{Toh&Yun2010}, \cite{Lee&Bresler2010}, \cite{RechtEtAl2010} and \cite{FornasierEtAl2011}.

For the general underdetermined setup, we first consider the aspect of designing an algorithm that can deal with any low-rank matrix. In addition, we consider the aspect of exploiting a-priori knowledge of underlying matrix structure for further performance improvement. The motivation for using structured low-rank matrices is due to their frequent occurrence in various signal processing and system identification problems.
We develop a new algorithm that addresses both aspects. Extending the approach of \cite{Tang&Nehorai2011}, the new algorithm is developed in an iterative framework based on least-squares (LS) estimation. For investigating matrix structure, we consider linearly structured matrices, such as Hankel and Toeplitz, as well as positive semidefinite matrices. The performance of the algorithm is evaluated by simulations and then compared to the Cramér-Rao bounds (CRBs).
The CRBs are presented for measurements in Gaussian noise. Specifically we derived the CRBs for Hankel and positive semidefinite matrices.

%

\emph{Notation:} The invertible vectorization and matrix construction
mappings are denoted $\text{vec}(\cdot): \mathbb{C}^{n \times p}
\rightarrow \mathbb{C}^{np
  \times 1}$ and $\text{mat}_{n,p}(\cdot): \mathbb{C}^{np \times 1}
\rightarrow \mathbb{C}^{n \times p}$, respectively. $\mathcal{X}_r \triangleq \{  \mathbf{X} \in \mathbb{C}^{n   \times p}: \text{rank}(\mathbf{X}) = r \}$ and $\mathcal{X}_+ \triangleq \{  \mathbf{X} \in \mathbb{C}^{n \times n}: \mathbf{X} \succeq \mathbf{0} \}$ denote the sets of rank $r$ matrices and positive semidefinite (p.s.d.) matrices, respectively.
Similarly, $\mathcal{X}_{S} \triangleq \{  \mathbf{X} \in
\mathbb{C}^{n   \times p}: \text{vec}(\mathbf{X}) = \mathbf{S}
\boldsymbol{\theta}, \boldsymbol{\theta} \in \mathbb{C}^{q}  \}$
denotes a subspace of linearly structured matrices specified by
$\mathbf{S} \in \mathbb{C}^{np \times q}$, which includes Hankel and Toeplitz
matrices. $\mathbf{A} \otimes \mathbf{B}$ is the Kronecker product
and $\langle \mathbf{A}, \mathbf{B} \rangle \triangleq \text{tr}
\left( \mathbf{B}^* \mathbf{A} \right)$ is the inner product.
$\mathbf{A}^\top$, $\mathbf{A}^*$ and $\mathbf{A}^\dagger$ denote the transpose, Hermitian transpose and Moore-Penrose pseudoinverse of $\mathbf{A} \in \mathbb{C}^{n
\times p}$, respectively while $\| \mathbf{A} \|_F$ is its Frobenius
norm. $\| \mathbf{x} \|_{\mathbf{W}} \triangleq \sqrt{ \mathbf{x}^* \mathbf{W} \mathbf{x} }$ is the weighted norm.

\section{System model}

\subsection{General underdetermined setup}

We consider a matrix $\mathbf{X} \in \mathcal{X}_r$, where $r \ll
\text{min}(n,p)$ is assumed to be known. It is observed by an undersampled linear
measurement process,
\begin{equation}
\mathbf{y} = \mathcal{A}(\mathbf{X}) + \mathbf{n} \in \mathbb{C}^{m
  \times 1},
\label{eq:measurement}
\end{equation}
where $m < np$ and the linear sensing operator $\mathcal{A}:
\mathbb{C}^{n \times  p} \rightarrow \mathbb{C}^{m \times 1}$ can be
written equivalently in forms,
\begin{equation}
\begin{split}
 \mathcal{A}( \mathbf{X} ) &=
 \begin{bmatrix}
 \langle \mathbf{X}, \mathbf{A}_1 \rangle \\
 \vdots \\
 \langle \mathbf{X}, \mathbf{A}_m \rangle
 \end{bmatrix}
= \mathbf{A}  \text{vec}(\mathbf{X}).
\end{split}
\label{eq:sensingoperator}
\end{equation}
The matrix $\mathbf{A}$ is assumed to be known and the measurement
noise $\mathbf{n}$ is assumed to be zero-mean, with $\E[
\mathbf{n}\mathbf{n}^*] = \mathbf{C} \in  \mathbb{C}^{m \times m}$ given.
Note that in matrix completion, $\left\{ \mathbf{A}_k \right\}$ is nothing but the set of element-selecting operators.

\subsection{Applications in signal processing and system identification}

Applications of the general underdetermined setup are illustrated in the
following examples:
\begin{enumerate}
\item  Recovery of data matrices $\mathbf{X} \in
\mathcal{X}_r$, compressed by some randomly chosen linear operation
$\mathcal{A}$ into $\mathbf{y}$.

\item Reconstruction of covariance matrices
$\text{Cov}(\mathbf{z}) = \mathbf{R} \in \mathcal{X}_r \cap
\mathcal{X}_+$ from a subset of second-order moments $r_{ij} = \E[
z_i z^*_j]$, estimated with zero-mean error $\varepsilon_{ij}$,
$(i,j) \in \Omega$. Then \eqref{eq:measurement} can be applied,
$\mathbf{y} = \mathcal{A}(\mathbf{R}) + \boldsymbol{\varepsilon}$.
In certain applications $\mathbf{R} \in \mathcal{X}_r \cap
\mathcal{X}_+ \cap \mathcal{X}_S$, e.g. Toeplitz or Persymmetric
structure.

\item Recovery of distance matrices $\mathbf{D} \in \mathcal{X}_r$, where $d_{ji} = d_{ij} = \| \mathbf{x}_i - \mathbf{x}_j \|^2_{\mathbf{W}}$ and $\mathbf{W} \succ \mathbf{0}$. A subset of distance measurements are observed in noise.

\item Identification of a system matrix $\mathbf{D} \in
  \mathcal{X}_r$. The system $\mathbf{z}_t = \mathbf{D} \mathbf{u}_t \in \mathbb{C}^{n
  \times 1}$ is sampled by a varying linear operator $\mathbf{y}_t =
\mathbf{A}_t \mathbf{z}_t + \mathbf{n}_t$, where $\mathbf{A}_t \in
\mathbb{C}^{s \times p}$. A special case is $\mathbf{A}_t =
\mathbf{e}^\top_{\ell(t)} $, where the index $\ell(t)$ varies the
sampling of $\mathbf{z}_t$ at each $t$. Given a set of input and
output data, $\{\mathbf{y}_t , \mathbf{u}_t  \}^T_{t=1}$, where $sT <
np$, the observations can be stacked by vectorization,
\begin{equation*}
\mathbf{y} =
\begin{bmatrix}
\mathbf{u}^\top_1  \otimes \mathbf{A}_1  \\
\vdots \\
\mathbf{u}^\top_T  \otimes \mathbf{A}_T
\end{bmatrix}
\text{vec}(\mathbf{D})  +  \mathbf{n}.
\end{equation*}
In certain scenarios $\mathbf{D}$ may also have linear structure.

\end{enumerate}

\section{Alternating least-squares estimator}

\subsection{General low-rank matrix reconstruction}

We begin by considering the case of reconstructing $\mathbf{X} \in
\mathcal{X}_r$. Using a weighted least-squares criterion, the
estimator is
\begin{equation}
\hat{\mathbf{X}} = \argmin_{\mathbf{X} \in \mathcal{X}_{r}} \|
\mathbf{y} - \mathcal{A}( \mathbf{X} ) \|^2_{\mathbf{C}^{-1}}.
\label{eq:ls}
\end{equation}
When the measurement noise $\mathbf{n}$ is Gaussian, this
estimator coincides with the maximum likelihood estimator. For
brevity we assume spatially uncorrelated noise, $\mathbf{C} =
\sigma^2 \mathbf{I}_m$, without loss of generality. Then minimizing the $\ell_2$-norm is equivalent to \eqref{eq:ls}. For general
$\mathbf{C}$ the observation model is pre-whitened by forming
$\bar{\mathbf{y}} = \mathbf{C}^{-1/2}\mathbf{y}$ and
$\bar{\mathbf{A}} = \mathbf{C}^{-1/2} \mathbf{A}$.

Since $\mathbf{X} \in \mathcal{X}_r$, we express $\mathbf{X} =
\mathbf{L}\mathbf{R}$ where $\mathbf{L} \in \mathbb{C}^{n \times
r}$ and $\mathbf{R} \in \mathbb{C}^{r \times p}$. Then the square of the measurement residual can be written as
\begin{equation}
\begin{split}
J(\mathbf{L},\mathbf{R}) &\triangleq \| \mathbf{y} - \mathcal{A}( \mathbf{L}\mathbf{R} ) \|^2_2 \\
&= \| \mathbf{y} - \mathbf{A}( \mathbf{I}_p \otimes \mathbf{L}
) \text{vec}(\mathbf{R})   \|^2_2 \\
&= \| \mathbf{y} - \mathbf{A}( \mathbf{R}^\top \otimes \mathbf{I}_n
) \text{vec}(\mathbf{L})   \|^2_2.
\end{split}
\end{equation}
The cost $J(\mathbf{L},\mathbf{R})$ is minimized in an alternating fashion by the following steps:
\begin{itemize}
 \item minimizing $\mathbf{R}$ with a fixed $\mathbf{L}$,
 \item minimizing $\mathbf{L}$ with a fixed $\mathbf{R}$.
\end{itemize}

In the new algorithm, the alternating minimization is performed through iterations. Starting with an initial $\mathbf{L}$, the iterations continue as long as the decreasing trend of $J(\mathbf{L},\mathbf{R})$ is observed.
Given $\mathbf{L}$, the minimizer of $\mathbf{R}$ is computed by $\text{vec}(\hat{\mathbf{R}}) = [ \mathbf{A} (
\mathbf{I}_p \otimes \mathbf{L} )  ]^\dagger \mathbf{y}$ and similarly, given $\mathbf{R}$, the minimizer of $\mathbf{L}$ is computed by $\text{vec}(\hat{\mathbf{L}})  = [\mathbf{A} ( \mathbf{R}^\top
\otimes \mathbf{I}_n )]^\dagger \mathbf{y}$. Since the Kronecker
products are sparse, $\mathbf{A} ( \mathbf{I}_p \otimes \mathbf{L}
) \in \mathbb{C}^{m \times rp}$ and $\mathbf{A} ( \mathbf{R}^\top
\otimes \mathbf{I}_n ) \in \mathbb{C}^{m \times nr}$ can be
computed efficiently. Henceforth, we refer to the algorithm as alternating LS (ALS).


\subsection{Structured low-rank matrix reconstruction}

Next, we consider a structured low-rank matrix $\mathbf{X} \in
\mathcal{X}_r$, and develop an ALS for a known matrix structure in Algorithm~\ref{alg:lse}.
In the algorithm, for each iteration, we approach the LS problem by first
relaxing the structural constraint, and compute $\mathbf{R}$ with a fixed
$\mathbf{L}$. Then, to impose the structural constraint on $\mathbf{R}$, the low-rank matrix estimate is projected onto the set of
structured matrices by $\bar{\mathbf{X}}  \triangleq \mathcal{P}(
\mathbf{L}\mathbf{R} )$, similar to `lift and project'
\cite{ChuEtAl2003}. $\mathbf{R}$ is subsequently modified as the
least-squares solution of $\bar{\mathbf{X}}$,
\begin{equation*}
\min_{\mathbf{R}} \| \mathbf{L} \mathbf{R} -
\bar{\mathbf{X}} \|^2_F.
\end{equation*}
$\mathbf{L}$ is updated in the same fashion. Here we mention that the algorithm is initialized by $\mathbf{L} := \mathbf{U}_0 \boldsymbol{\Sigma}_0$ where $[\mathbf{U}_0, \boldsymbol{\Sigma}_0, \mathbf{V}_0]$ =
 \texttt{svdtrunc}$ \left( \text{mat}_{n,p}(\mathbf{A}^* \mathbf{y}), r  \right)$ and \texttt{svdtrunc}$\left( \mathbf{Z}, r
\right)$ denotes the singular value decomposition of
$\mathbf{Z} \in \mathbb{C}^{n \times p}$ truncated to the $r$th
singular value.


\begin{algorithm}
\caption{: ALS with known matrix structure} \label{alg:lse}
\begin{algorithmic}[1]
 \State  Input: $\mathbf{y}$, $\mathbf{A}$ and $r$
 \State $[\mathbf{U}_0, \boldsymbol{\Sigma}_0, \mathbf{V}_0]$ =
 \texttt{svdtrunc}$ \left( \text{mat}_{n,p}(\mathbf{A}^* \mathbf{y}), r  \right)$
 \State $\mathbf{L} := \mathbf{U}_0 \boldsymbol{\Sigma}_0$
 \While{$J(\mathbf{L},\mathbf{R})$ decreases}
 \State $\mathbf{R} := \text{mat}_{r,p} \left(   [ \mathbf{A} (
\mathbf{I}_p \otimes \mathbf{L} )  ]^\dagger \mathbf{y} \right)$
  \State $\bar{\mathbf{X}} := \mathcal{P}(\mathbf{L} \mathbf{R} )$
 \State $\mathbf{R} := \mathbf{L}^\dagger \bar{\mathbf{X}} $
  \State $\mathbf{L} := \text{mat}_{n,r} \left(  [\mathbf{A} (
   \mathbf{R}^\top \otimes \mathbf{I}_n )]^\dagger \mathbf{y}
 \right)$
  \State $\bar{\mathbf{X}} := \mathcal{P}(\mathbf{L} \mathbf{R} )$
 \State $\mathbf{L} := \bar{\mathbf{X}} \mathbf{R}^\dagger$
 \EndWhile
 \State Output: $\hat{\mathbf{X}} = \mathbf{LR}$
\end{algorithmic}
\end{algorithm}

For linearly structured matrices the projection, $\bar{\mathbf{X}}
= \mathcal{P}_{\mathcal{X}_{S}}(\mathbf{LR})$, is computed by
obtaining $\bar{\boldsymbol{\theta}} = \mathbf{S}^\dagger
\text{vec}(\mathbf{LR})$ and setting $\bar{\mathbf{X}} =
\text{mat}_{n,p}(\mathbf{S} \bar{\boldsymbol{\theta}})$. For
p.s.d. matrices, $\bar{\mathbf{X}} =
\mathcal{P}_{\mathcal{X}_{+}}(\mathbf{LR})$ is computed by first
performing an eigenvalue decomposition of a symmetrized matrix,
which projects the matrix onto the set of Hermitian matrices,
$\frac{1}{2} ( \mathbf{LR} + (\mathbf{LR})^*) = \mathbf{V}
\boldsymbol{\Lambda} \mathbf{V}^*$. The decomposition is performed to the $r$th largest eigenvalue, then setting $\bar{\mathbf{X}}
= \mathbf{V}_r \widetilde{\boldsymbol{\Lambda}}_r \mathbf{V}^*_r$, where positive eigenvalues have been retained, projects the matrix onto $\mathcal{X}_+$.

\section{Cramér-Rao bounds}

For relevant comparison, we use CRBs. In this section, we describe the CRB expressions.

\subsection{Unstructured matrix}

For real-valued $\mathbf{X}$ and Gaussian measurement noise
$\mathbf{n}$, $\mathbf{y}$ is distributed as $\mathcal{N}( \mathbf{A}\text{vec}(\mathbf{X}), \mathbf{C} )$.
The Cramér-Rao bound for unbiased low-rank matrix estimators
without structure was derived in \cite{Tang&Nehorai2011} (see also \cite{Werner&Jansson2006}):
$\E_{y|X}[ \| \mathbf{X} - \hat{\mathbf{X}}(\mathbf{y}) \|^2_F ]
\geq \text{CRB}(\mathbf{X})$, where
\begin{equation}
\text{CRB}(\mathbf{X}) = \text{tr}\left(  (\mathbf{P}^\top
\mathbf{A}^\top\mathbf{C}^{-1}\mathbf{A}\mathbf{P} )^{-1} \right).
\label{eq:CRLB_Nehorai}
\end{equation}
The CRB holds when $\text{rank}(\mathbf{AP}) = r(n+p)-r^2$ where
$\mathbf{P} \triangleq \begin{bmatrix} \mathbf{V}_1 \otimes
  \mathbf{U}_0  & \mathbf{U}_0 \otimes \mathbf{V}_0 & \mathbf{V}_0 \otimes
  \mathbf{U}_1 \end{bmatrix}$.
The submatrices are obtained from the
left-singular vectors, $\mathbf{U}_0 = [\mathbf{u}_1, \dots,
\mathbf{u}_r]$ and $\mathbf{U}_1 = [\mathbf{u}_{r+1}, \dots,
\mathbf{u}_{n}]$, and right-singular vectors, $\mathbf{V}_0 =
[\mathbf{v}_1, \dots, \mathbf{v}_r]$ and $\mathbf{V}_1 =
[\mathbf{v}_{r+1}, \dots, \mathbf{v}_{n}]$, of $\mathbf{X} \in \mathcal{X}_r$.

\subsection{Structured matrices}

\subsubsection{Hankel and Toeplitz matrices}
For certain linearly structured matrices, such as Hankel or Toeplitz, low-rank parametrizations $\mathbf{S} \boldsymbol{\theta} = \mathbf{S}\mathbf{g}(\boldsymbol{\alpha})$ exist that enable the derivation of CRBs. E.g. when $\mathcal{X}_S$ is Hankel, $\mathbf{X}$ can be parameterized in the \emph{canonical controllable form}
\cite{Kailath1980} as $\left[ \mathbf{X} \right]_{ij} =
\mathbf{b}^\top \boldsymbol{\Phi}^{i+j-2} \mathbf{e}_1$, where
$\mathbf{b} \in \mathbb{R}^{r}$,
\begin{equation*}
\begin{split}
\boldsymbol{\Phi} =
\begin{bmatrix}
-a_1 & -a_2 & \dots & -a_{r-1} & -a_r\\
1 & 0 & \dots & 0 & 0\\
0 & 1 & \dots & 0 & 0\\
\vdots & & \ddots & & \vdots\\
0 & 0 & \dots & 1 & 0
\end{bmatrix}
\end{split} \in \mathbb{R}^{r \times r},
\end{equation*}
and $\mathbf{e}_1$ is the first standard basis vector in $\mathbb{R}^r$.
Here the parameters are $\boldsymbol{\alpha} = \begin{bmatrix}
\mathbf{a}^\top & \mathbf{b}^\top \end{bmatrix}^\top$, where
$\mathbf{a} = \begin{bmatrix} a_1 & \dots a_r \end{bmatrix}^\top$. A similar parametrization can be found for Toeplitz matrices.

The Cramér-Rao bound is then given by
\begin{equation}
\text{CRB}_S(\mathbf{X}) = \text{tr}\left(
\boldsymbol{\Delta}_{\boldsymbol{\alpha}}^\top
\mathbf{J}^{-1}(\boldsymbol{\alpha})
\boldsymbol{\Delta}_{\boldsymbol{\alpha}} \right)
\label{eq:CRLB_linstruct}
\end{equation}
where $\boldsymbol{\Delta}_{\boldsymbol{\alpha}} = \frac{\partial
\mathbf{g}(\boldsymbol{\alpha})}{\partial \boldsymbol{\alpha}
}^\top \mathbf{S}^\top$ and $\mathbf{J}(\boldsymbol{\alpha})
\triangleq - \E \left[ \nabla^2_{\boldsymbol{\alpha}}
p(\mathbf{y}; \alpha)
 \right]$ is the Fisher information
matrix of $\boldsymbol{\alpha}$ \cite{Kay1993},
\begin{equation*}
\mathbf{J}(\boldsymbol{\alpha}) =
\boldsymbol{\Delta}_{\boldsymbol{\alpha}} \mathbf{A}^\top
\mathbf{C}^{-1} \mathbf{A}
\boldsymbol{\Delta}^\top_{\boldsymbol{\alpha}}.
\end{equation*}

\subsubsection{Positive semidefinite matrix}
Positive semidefinite matrices can be parameterized as $\mathbf{X}
= \mathbf{M}\mathbf{M}^\top$, where $\mathbf{M} \in \mathbb{R}^{n
\times r}$. Let $\boldsymbol{\alpha} = \text{vec}(\mathbf{M})$, so
that $\mathbf{X} = \mathbf{g}(\boldsymbol{\alpha})$ then
$\mathbf{J}(\boldsymbol{\alpha})$ has the same form as above, with  $\boldsymbol{\Delta}_{\boldsymbol{\alpha}} = \frac{\partial
\mathbf{g}(\boldsymbol{\alpha})}{\partial \boldsymbol{\alpha}
}^\top $. The gradient can be written compactly as $\frac{\partial
\mathbf{g}(\boldsymbol{\alpha})}{\partial \boldsymbol{\alpha} } =
\left[ (\mathbf{I}_n \otimes \mathbf{M})\mathbf{T} + (\mathbf{M}
\otimes \mathbf{I}_n ) \right]$ where $\mathbf{M} =
\text{mat}_{n,r}(\boldsymbol{\alpha})$ and $\mathbf{T}$ is the
matrix form of the transpose operator, $\mathbf{T}
\text{vec}(\mathbf{M}) = \text{vec}(\mathbf{M}^\top)$. Note that
parametrization by $\mathbf{M}$ is unique only up to an
orthonormal transformation, hence
$\mathbf{J}(\boldsymbol{\alpha})$ is not invertible in general.
The Cramér-Rao bound is then given by
\begin{equation}
\text{CRB}_+(\mathbf{X}) = \text{tr}\left(
\boldsymbol{\Delta}^\top_{\boldsymbol{\alpha}}
\mathbf{J}^{\dagger}(\boldsymbol{\alpha})
\boldsymbol{\Delta}_{\boldsymbol{\alpha}} \right)
\label{eq:CRLB_posdef}
\end{equation}
provided that
\begin{equation*}
\boldsymbol{\Delta}^\top_{\boldsymbol{\alpha}} =
\boldsymbol{\Delta}^\top_{\boldsymbol{\alpha}}
\mathbf{J}(\boldsymbol{\alpha})
\mathbf{J}^\dagger(\boldsymbol{\alpha})
\end{equation*}
holds, or equivalently that $\mathbf{P}^\perp_{\mathbf{J}}
\boldsymbol{\Delta}_{\boldsymbol{\alpha}} = \mathbf{0}$, where
$\mathbf{P}^\perp_{\mathbf{J}}$ is the orthogonal projection onto
the nullspace of $\mathbf{J}(\boldsymbol{\alpha})$
\cite{Werner&Jansson2006}.

For sake of brevity the CRB derivations for Hankel and p.s.d. matrices are given in a supplementary note \cite{ZachariahEtAl2011}.

\section{Experiments and Results}


\subsection{Experiment setup and performance measure}

In the following we consider real-valued matrices of dimension $n = p \equiv 100$, with Hankel and p.s.d. structure respectively. For Hankel structure we generate $\mathbf{X}$ randomly by first creating a matrix with elements from an i.i.d. $\mathcal{N}(0,1)$ and fitting
$\boldsymbol{\alpha}$ using Prony's method \cite{Hayes1996}. Then set $\mathbf{X} =
\text{mat}_{n,p}( \mathbf{S} \mathbf{g}(\boldsymbol{\alpha}) )$. For
p.s.d. structure, $\mathbf{X}$ is generated by $\mathbf{X} = \mathbf{M} \mathbf{M}^\top$, where the elements of $\mathbf{M}$ are generated by  i.i.d. $\mathcal{N}(0,1)$.
We let parameter $\lambda \triangleq r / \min(n,p) \in (0,1]$ determine the rank, which controls the degrees of freedom.

The linear sensing operator, modeled by $\mathbf{A}$
in \eqref{eq:sensingoperator}, is selected randomly by $a_{ij}
\sim \mathcal{N}(0,\frac{1}{m})$ \cite{RechtEtAl2010}. The
measurement noise is generated as $\mathbf{n} \sim
\mathcal{N}(\mathbf{0}, \mathbf{\sigma}^2 \mathbf{I}_m)$. In the experiments, two signal parameters are varied:
\begin{enumerate}

\item Signal to measurement noise ratio,
\begin{equation}
\text{SMNR} \triangleq \frac{\E \left[ \| \mathbf{X} \|^2_F
\right] }{\E \left[ \| \mathbf{n} \|^2_F \right] }.
\end{equation}

\item Measurement fraction, $\rho \in (0,1]$, so that $m \equiv
\lceil \rho np \rceil$.

\end{enumerate}

The signal-to-reconstruction error ratio, or inverse of the normalized mean square error (NMSE),
\begin{equation}
\text{SRER} \triangleq \frac{\E \left[ \| \mathbf{X} \|^2_F
  \right] }{\E \left[ \| \mathbf{X} - \hat{\mathbf{X}}\|^2_F \right]
  } = \frac{1}{\text{NMSE}},
\end{equation}
is used as the performance measure as it increases with SMNR. For an unbiased estimator, $\text{SRER}  \leq \E [ \| \mathbf{X} \|^2_F  ] / \E [ \text{CRB}(\mathbf{X}) ] $. When SRER = 0~dB there is no reconstruction gain over the zero solution $\hat{\mathbf{X}} =
\mathbf{0}$.

\subsection{Results}

The experiments were repeated for 500 Monte Carlo simulations. For each run a new realization of $\mathbf{X}$, $\mathbf{y}$ and $\mathbf{A}$ was generated and the Cramér-Rao bounds were computed correspondingly. The algorithm was set to terminate when the measurement residual $J(\mathbf{L},\mathbf{R})$ stops decreasing.

\begin{figure*}
  \begin{center}
    \includegraphics[width=1.97\columnwidth]{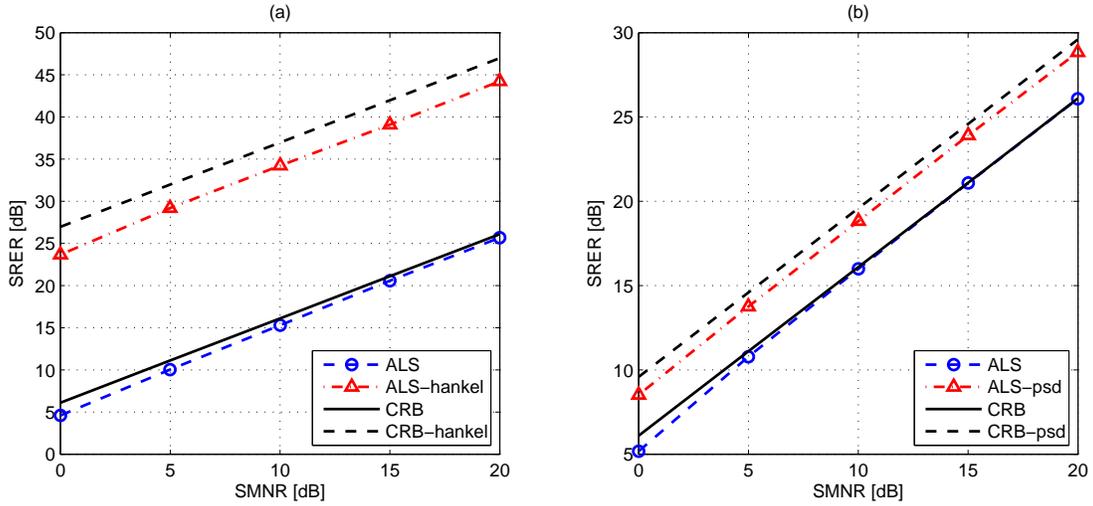}
  \end{center}
  \vspace*{-0.5cm}
  \caption{SMNR versus SRER for $n=p=100$, $\lambda = 0.03$ and $\rho = 0.3$. (a) Hankel structure, (b) Positive semidefinite structure.}
  \label{fig:smnr}
\end{figure*}

Figure \ref{fig:smnr} shows the performance of ALS when varying SMNR for Hankel and p.s.d. matrices, respectively. The measurement fraction was fixed to $\rho = 0.3$ and
the relative rank was $\lambda = 0.03$. The algorithm is tested both with and without prior knowledge of matrix structure. Without such information it quickly approaches the bound for unstructured matrices \eqref{eq:CRLB_Nehorai}. The reconstruction gain is significantly raised when exploiting the matrix structure. ALS remains at a SRER gap from the bounds \eqref{eq:CRLB_linstruct} and \eqref{eq:CRLB_posdef} because it relaxes the problem by alternating projections. The gaps are 2.75 and 0.77~dB for Hankel and p.s.d., respectively, at SMNR = 10~dB. The estimator inserts a proportionally larger bias into the MSE compared to the estimator without prior information.

\begin{figure*}
  \begin{center}
    \includegraphics[width=1.97\columnwidth]{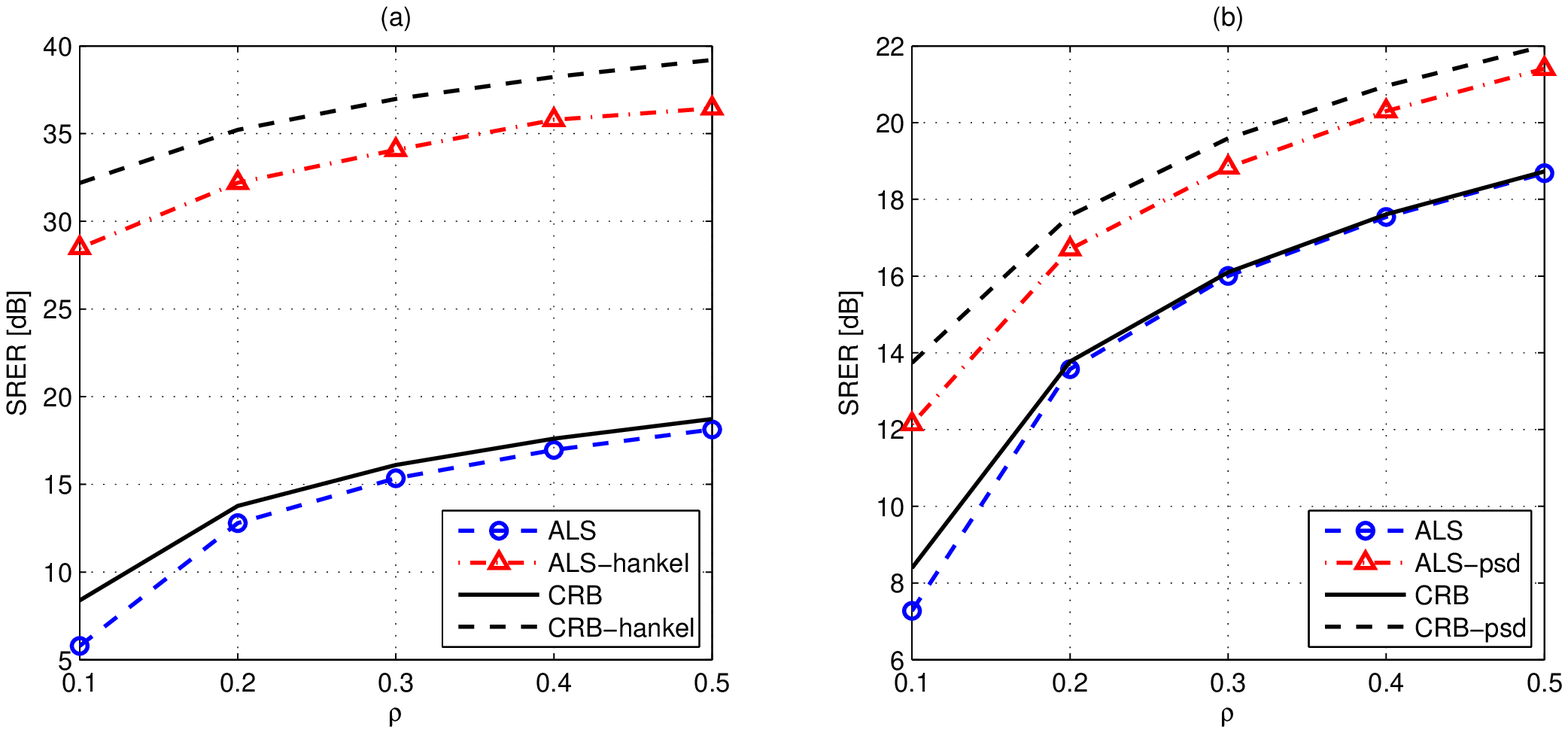}
  \end{center}
  \vspace*{-0.5cm}
  \caption{$\rho$ versus SRER for $n=p=100$, $\lambda = 0.03$ and $\text{SMNR} = 10$~dB. (a) Hankel structure, (b) Positive semidefinite structure.}
  \label{fig:rho}
\end{figure*}

In Figure \ref{fig:rho} the experiment is repeated for varying $\rho$, while $\lambda = 0.03$ and SMNR is fixed to 10~dB. As the measurement fraction increases the SRER of the algorithm without prior matrix structure approaches the corresponding CRB. At very low $\rho$, below 0.1, the bound \eqref{eq:CRLB_Nehorai} tends to break down as the rank constraint is violated. Given prior structural information the algorithm performs similar to the SMNR case.

\emph{Reproducible research:} The \textsc{Matlab} code to reproduce the figures in this paper are available at http://sites.google.com/site/saikatchatt/, along with the figures for varying $\lambda$ that have been omitted here due to lack of space.

\vspace*{-0.25cm}

\section{Conclusions}

We have developed an algorithm based on least-squares estimation for reconstruction of low-rank matrices in a general underdetermined setup. Furthermore it is capable of exploiting structures, in particular linearly structured matrices and positive semidefinite matrices, leading to better performance. Through simulations, we found that the algorithm provides good performance compared to the Cramér-Rao bound.

\bibliography{refs_lowrank}
\bibliographystyle{ieeetr}

\end{document}